\documentclass[12pt, reqno]{amsart}
\allowdisplaybreaks
\begin{document}
\setcounter{page}{1}

\title[\hfilneg \hfil Iyengar's  inequality ]
{$\pmb{\psi}$-Caputo Fractional Iyengar's Type Inequalities}

\author[Bhagwat Yewale \hfil]
{Bhagwat R. Yewale}

\address{Bhagwat R. Yewale \newline
 Dept. of Mathematics,
 Dr. B. A. M. University, Aurangabad,
 Maharashtra 431004, India}
\email{yewale.bhagwat@gmail.com}

\author[Deepak Pachpatte \hfil]
{Deepak B. Pachpatte}

\address{Deepak B. Pachpatte \newline
 Dept. of Mathematics,
 Dr. B. A. M. University, Aurangabad,
 Maharashtra 431004, India}
\email{pachpatte@gmail.com}

\subjclass[2010]{26A33, 26D10, 26D15}
\keywords{Iyengar Inequality, $\psi$-Caputo fractional derivative, Taylor's formulae}

\begin{abstract}
  In this paper, we establish Iyengar type inequalities utilizing $\psi$-Caputo fractional derivatives that is, fractional derivative of a function with respect to another function, which is generalization of some known fractional derivatives such as Riemann-Liouville, Hadamard, Erd\'{e}lyi-Kober. The inequalities in this article are with respect to $L_{p}$ norms, $1\leq p \leq \infty.$ The tools used in the analysis are based on Taylor's formulae for $\psi$-Caputo fractional derivatives.

\end{abstract}

\maketitle

\section{Introduction}
In the past several years, the development of the theories of fractional calculus is well contributed by many researchers. The field of fractional calculus mainly focuses on the study of integrals and derivatives of non integer order. Also, this field has been attracting many researchers due to the scope of research as well as providing significant applications in various fields, for example control theory \cite{Pet}, bioengineering \cite{Mag}, viscoelasticity \cite{Koe}. Development of the theory of fractional calculus leads to the various fractional operators (i.e. fractional integrals and fractional derivative), such as the Riemann-Liouville, Caputo, Hadamard, Erd\'{e}lyi-Kober, Riesz (see \cite{Kil}).

In literature, there are numerous approaches acquiring a generalization of fractional operators that authorize the exhibition of distinct kernel. For instance, Prabhakar \cite{Pra}, introduced a new fractional integral operator with the help of three parameter Mittag-Leffler function. In \cite{Sri}, Srivastava with Tomovski investigated generalized Mittag-Leffler function and defined corresponding fractional integral operator known as the generalized Prabhakar fractional integral. Almeida \cite{Alm1} have introduced one more dimension to the study of fractional calculus by presenting a fractional derivative of a function with respect to another function in Caputo sense known as $\psi$-Caputo fractional derivative. In \cite{Alm2, Alm3, Awa} authors provides applications and various results for $\psi$-Caputo fractional derivative. Recently, motivated by the definition of Hilfer fractional dervative Sousa and Oliveira \cite{Sou}, established new fractional derivative called $\psi$-Hilfer fractional derivative.

A wide variety of such operators led researchers to use them, however some authors studied fractional differential problems and some authors established fractional inequalities for various fractional operators. For instance, Pachpatte \cite{Pac1}, studied \v{C}eby\v{s}ev like inequalities for functions of two and three variables using $\psi$-Caputo fractional derivative. In \cite{Esh}, Eshaghi with Ansari established Lyapunov inequality for fractional differential equations for Prabhakar derivative. Rashid et al. \cite{Ras1, Ras2} investigated new fractional integral operator known Generalized proportional fractional integral operators with respect to another function and established certain Gr\"{u}ss type and the reverse Minkowski type inequalities for this operators. Aljaaidi with Pachpatte \cite{Alj}, presented Minkowski inequalities by mean of $\psi$-Riemann Liouville fractional integral operators. A large variety of inequalities have been proposed and studied for different fractional operators, see \cite{Bud, Pur, Sar2, Set} and the references therein. In last few decades, a number of mathematician have devoted their efforts to generalize, extends and give numerous variants of the inequalities associated with the names of Minkowski, H\"{o}lder's, Hardy, Trapezoid, Ostrowski, \v{C}eby\v{s}ev, for more details, see \cite{Mit1, Pac2}. In 1938, Iyengar proved the following useful inequality (see \cite{Mit2}, p. 471):

Let $\mathrm{\textsf{g}}$ be a differentiable function on [a, b]. If $|\mathrm{\textsf{g}}'(z)| \leq \mathrm{M},\forall z\in$ [a, b], then following inequality holds:
\[\label{eq1.1}
\bigg|\int_{a}^{b}\mathrm{\textsf{g}}(z)dz-\frac{1}{2}(b-a)(\mathrm{\textsf{g}}(a)+\mathrm{\textsf{g}}(b))\bigg|
\leq\frac{\mathrm{M}(b-a)^{2}}{4}-
\frac{(\mathrm{\textsf{g}}(b)-\mathrm{\textsf{g}}(a))^{2}}{4\mathrm{M}}.\tag{1.1}
\]

Due to its usefulness, over the past four decades, inequality \eqref{eq1.1} has been extensively studied and generalized in various directions as well as some of its applications can be found in \cite{Aga2, Che, Mil, Vas}. Agarwal with Dragomir \cite{Aga1}, presented an extension of the inequality \eqref{eq1.1} by using Steffensen's inequality as a main tool. However it has been seen that the Iyengar type inequalities established by Elezovi\'{c} with Pe\v{c}ari\'{c} \cite{Ele} and \v{C}uljak with Elezovi\'{c} \cite{Cul} had used the same tool and also provided more improvements in it. Furthermore, in \cite{Liu1}, Liu obtained the Iyengar type inequalities for certain weaker conditions that is, the function in the Iyengar type inequality is not necessarily differentiable. Anastassiou \cite{Ana1,Ana2}, studied the above mentioned inequality by employing Caputo fractional derivative and Canavati fractional derivative. In recent study, significant contribution have been made to which inequality \eqref{eq1.1} leads, and numerous extensions and improvements can be found related to it in the literature \cite{Ana3, Dra, Fra1, Fra2, Liu2, Sar1}.

Motivated by the aforementioned work, in this paper, we derive the Iyengar type inequalities involving Caputo operator with respect to the new function $\psi$. The overall structure of this paper take the form of three sections with an introduction. The remaining paper is organized as follows: In section 2, we present some preliminaries, essential definitions and mathematical tools which are used to carry out our work. Main results of the Iyengar type inequalities are established in the section 3.

\section{Preliminaries}
In this section we give some basic definitions, preliminaries and weighted spaces which are useful for our subsequent discussions.
Let [a, b] with ($-\infty<a<b<+\infty$) be a finite interval in $\mathbb{R^{+}}$.  Let $\mathcal{\textit{A}C}^{n}([a, b])$ be the space of n-times absolutely continuous differentiable functions on [a, b], $\mathcal{C}^{n}$ be the space of n-times continuously differentiable functions on [a, b] and $\mathcal{L}_{p}([a, b])$, $1\leq p \leq \infty$ be the space of Lebesgue integrable functions endowed with the norm
\[
\|\mathrm{\textsf{h}}\|_{\mathcal{L}_{p}([a, b])}= \underset{\mathrm{x}\in[a, b]}{\mathrm{max}}{|\mathrm{\textsf{h}(x)}|}.
\]
The space of n-times absolutely continuous differentiable function $\mathrm{\textsf{h}}$ on [a, b] is defined as\newline
\[
\mathcal{\textit{A}C}^{n}([a, b])=\{\mathrm{\textsf{h}}:[a, b]\rightarrow\mathbb{R}:\mathrm{\textsf{h}}^{(n-1)}
\in\mathcal{\textit{A}C}([a, b])\}.
\]

Now we highlight some definitions and Taylor's formulae of a fractional integrals and derivative of a function with respect to another function that will be useful for the development of our article:\newline
{\textbf{Definition 2.1}}. \cite{Kil} Let $\alpha>0,$ \textit{f} be an integrable function on [a, b] and $\psi \in\mathcal{C}^{1}[a, b]$ be an increasing function such that $\psi^{'}(t)\neq 0,$ $ \forall t \in [a, b]$ then left sided $\psi$-Riemann-Liouvlle fractional integral of a function $\textit{f}$ is given by
\[
\mathcal{I}^{\alpha, \psi}_{a+}f(t)=\frac{1}{\Gamma(\alpha)}\int_{a}^{t}\psi^{'}(s)(\psi(t)-\psi(s))^{\alpha-1}f(s)ds.
\]
{\textbf{Definition 2.2}}. \cite{Kil} Let $\alpha>0,$ \textit{f} be an integrable function on [a, b] and $\psi \in\mathcal{C}^{1}[a, b]$ be an increasing function such that $\psi^{'}(x)\neq 0,$ $\forall$ $ t \in [a, b]$ then right sided $\psi$-Riemann-Liouville fractional integral of a function $\textit{f}$ is given by
\[
\mathcal{I}^{\alpha, \psi}_{b-}f(t)=\frac{1}{\Gamma(\alpha)}\int_{t}^{b}\psi^{'}(s)(\psi(s)-\psi(t))^{\alpha-1}f(s)ds.
\]
{\textbf{Definition 2.3}}. \cite{Alm1} Let $\alpha>0, n\in \mathbb{N}$. Let \textit{f}, $\psi \in\mathcal{C}^{n}[a, b]$ with $\psi$ an increasing function such that $\psi^{'}(t)\neq 0,$ $\forall$ $t \in [a, b]$ then left sided $\psi$-Caputo fractional derivative of a
function $\textit{f}$ is defined as
\[
\mathcal{D}^{\alpha,\psi}_{a+}f(t)=\frac{1}{\Gamma(n-\alpha)}\int_{a}^{t}\psi^{'}(s)(\psi(t)-\psi(s))^{n-\alpha-1}
\bigg(\frac{1}{\psi^{'}(t)}\frac{d}{dt}\bigg)^{n}f(s)ds.
\]
Where $n=\lceil \alpha \rceil+1$, $\lceil \alpha \rceil$ is the integer part of $\alpha$.\newline
{\textbf{Definition 2.4}}. \cite{Alm1} Let $\alpha>0, n\in \mathbb{N}$. Let \textit{f}, $\psi \in\mathcal{C}^{n}[a, b]$ with $\psi$ an increasing function such that $\psi^{'}(t)\neq 0,$ $\forall$ $t \in [a, b]$ then right sided $\psi$-Caputo fractional derivative of a
function $\textit{f}$ is defined as
\[
\mathcal{D}^{\alpha,\psi}_{b-}f(t)=\frac{1}{\Gamma(n-\alpha)}\int_{t}^{b}\psi^{'}(s)(\psi(s)-\psi(t))^{n-\alpha-1}
\bigg(-\frac{1}{\psi^{'}(t)}\frac{d}{dt}\bigg)^{n}f(s)ds.
\]
Where $n=\lceil \alpha \rceil+1$, $\lceil \alpha \rceil$ is the integer part of $\alpha$.\newline
{\textbf{Lemma 2.1}}. \cite{Alm1} The left sided and right sided $\psi$-fractional Taylor's formulae are given as follows:
\[
f(t)=\sum_{k=0}^{n-1}\frac{f^{[k]}_{\psi}(a)}{k!}(\psi(t)-\psi(a))^{k}
+\mathcal{I}^{\alpha, \psi}_{a+}\mathcal{D}^{\alpha,\psi}_{a+}f(t)
\]
and
\[
f(t)=\sum_{k=0}^{n-1}(-1)^{k}\frac{f^{[k]}_{\psi}(b)}{k!}(\psi(b)-\psi(t))^{k}
+\mathcal{I}^{\alpha, \psi}_{b-}\mathcal{D}^{\alpha,\psi}_{b-}f(t),
\]
where $f^{[k]}_{\psi}(t)=\big(\frac{1}{\psi^{'}(t)} \frac{d}{dt}\big)^{k}f(t)$.
\section{Iyengar Inequalities using $\psi$-Caputo Fractional Operators}
This section deals with the Iyengar type inequalities via $\psi$-Caputo fractional derivative operators.\newline
{\textbf{Theorem 3.1}} Let $\alpha>0$, $\textit{f}\in \mathcal {\textit{A}C}^{n}([a, b])$ and $\psi \in \mathcal{C}^{n}([a, b])$ with $\psi$ is an increasing and $\psi^{'}(t)\neq0$, $\forall t\in [a, b]$. Suppose $\mathcal{D}^{\alpha, \psi}_{a+}f$, $\mathcal{D}^{\alpha, \psi}_{b-}f \in \mathcal{L}_{\infty}([a, b])$. Then for $a\leq s \leq b$ following inequalities hold:\newline
(i)
\begin{align*}\label{eq3.1}
&\bigg|\int_{a}^{b}f(t)dt-\sum_{k=0}^{n-1}\frac{1}{(k+1)!}
\big(f^{[k]}_{\psi}(a)(\psi(s)-\psi(a))^{k+1}+(-1)^{k}f^{[k]}_{\psi}(b)(\psi(b)-\psi(s))^{k+1}\big)\bigg| \\
\leq \quad& \frac{\mathrm{max}\big[\|\mathcal{D}^{\alpha,\psi}_{a+}f\|_{L_{\infty}([a, b])},
\|\mathcal{D}^{\alpha,\psi}_{b-}f\|_{L_{\infty}([a, b])}\big]}
{\Gamma(\alpha+2)}[(\psi(s)-\psi(a))^{\alpha+1}+(\psi(b)-\psi(s))^{\alpha+1}].\tag{3.1}
\end{align*}
(ii)
The right hand side of the inequality \eqref{eq3.1} is minimized at $\psi(s)=\frac{\psi(a)+\psi(b)}{2}$, with the value
$\frac{\mathrm{max}\big[\|\mathcal{D}^{\alpha,\psi}_{a+}f\|_{L_{\infty}([a, b])},
\|\mathcal{D}^{\alpha,\psi}_{b-}f\|_{L_{\infty}([a, b])}\big]}{\Gamma(\alpha+2)}\frac{(\psi(b)-\psi(a))^{\alpha+1}}{2^{\alpha}}$, that is
\begin{align*}\label{eq3.2}
\bigg|\int_{a}^{b}f(t)dt-&\sum_{k=0}^{n-1}\frac{1}{(k+1)!}\frac{(\psi(b)-\psi(a))^{k+1}}{2^{k+1}}
\big(f^{[k]}_{\psi}(a)+(-1)^{k}f^{[k]}_{\psi}(b)\big)\bigg|\\
\leq \quad& \frac{\mathrm{max}\big[\|\mathcal{D}^{\alpha,\psi}_{a+}f\|_{L_{\infty}([a, b])},
\|\mathcal{D}^{\alpha,\psi}_{b-}f\|_{L_{\infty}([a, b])}\big]}
{\Gamma(\alpha+2)}\frac{(\psi(b)-\psi(a))^{\alpha+1}}{2^{\alpha}}.\tag{3.2}
\end{align*}
(iii)
When $f^{[k]}_{\psi}(a)=f^{[k]}_{\psi}(b)=0$, for $k=0,1,2,...,n-1,$ from \eqref{eq3.2} we get a sharp inequality as follows:
\begin{align*}\label{eq3.3}
\bigg|\int_{a}^{b}f(t)dt\bigg|
\leq \quad& \frac{\mathrm{max}\big[\|\mathcal{D}^{\alpha,\psi}_{a+}f\|_{L_{\infty}([a, b])},
\|\mathcal{D}^{\alpha,\psi}_{b-}f\|_{L_{\infty}([a, b])}\big]}
{\Gamma(\alpha+2)}\frac{(\psi(b)-\psi(a))^{\alpha+1}}{2^{\alpha}}.\tag{3.3}
\end{align*}
(iv)
In general, for $i=0,1,2,...m\in\mathbb{N}$, we obtain
\begin{align*}\label{eq3.4}
\bigg|\int_{a}^{b}f(t)dt-\sum_{k=0}^{n-1}&\frac{1}{(k+1)!}\bigg(\frac{\psi(b)-\psi(a)}{m}\bigg)^{k+1}\big[i^{k+1}f^{[k]}_{\psi}(a)
+(-1)^{k}(m-i)^{k+1}f^{[k]}_{\psi}(b)\big]\bigg|\\
\leq&\quad \frac{\mathrm{max}\big[\|\mathcal{D}^{\alpha,\psi}_{a+}f\|_{L_{\infty}([a, b])},
\|\mathcal{D}^{\alpha,\psi}_{b-}f\|_{L_{\infty}([a, b])}\big]}
{\Gamma(\alpha+2)}\bigg(\frac{\psi(b)-\psi(a)}{m}\bigg)^{\alpha+1}\\
&\quad\big[i^{\alpha+1}+(m-i)^{\alpha+1}\big].\tag{3.4}
\end{align*}
(v)
If $f^{[k]}_{\psi}(a)=f^{[k]}_{\psi}(b)=0, k=1,2,...,n-1$, then from \eqref{eq3.4} we obtain
\begin{align*}\label{eq3.5}
\bigg|\int_{a}^{b}f(t)dt-&\bigg(\frac{\psi(b)-\psi(a)}{m}\bigg)\big[if(a)+(m-i)f(b)\big]\bigg|\\
\leq &\quad\frac{\mathrm{max}\big[\|\mathcal{D}^{\alpha,\psi}_{a+}f\|_{L_{\infty}([a, b])},
\|\mathcal{D}^{\alpha,\psi}_{b-}f\|_{L_{\infty}([a, b])}\big]}
{\Gamma(\alpha+2)}\bigg(\frac{\psi(b)-\psi(a)}{m}\bigg)^{\alpha+1}\\
&\quad \big[i^{\alpha+1}+(m-i)^{\alpha+1}\big].\tag{3.5}
\end{align*}
(vi)
For {\it i=1, m=2} from \eqref{eq3.5} we get
\begin{align*}\label{eq3.6}
\bigg|\int_{a}^{b}f(t)dt-&\bigg(\frac{\psi(b)-\psi(a)}{2}\bigg)\big[f(a)+f(b)\big]\bigg|\\
\leq &\quad \frac{\mathrm{max}\big[\|\mathcal{D}^{\alpha,\psi}_{a+}f\|_{L_{\infty}([a, b])},
\|\mathcal{D}^{\alpha,\psi}_{b-}f\|_{L_{\infty}([a, b])}\big]}
{\Gamma(\alpha+2)}\frac{\big(\psi(b)-\psi(a)\big)^{\alpha+1}}{2^{\alpha}}.\tag{3.6}
\end{align*}
{\textbf{Proof}} (i)\newline
From the left $\psi$-Caputo fractional Taylor's formula we have
\[\label{eq3.7}
f(t)-\sum_{k=0}^{n-1}\frac{f^{[k]}_{\psi}(a)}{k!}(\psi(t)-\psi(a))^{k}
=\frac{1}{\Gamma(\alpha)}\int_{a}^{t}\psi^{'}(s)(\psi(t)-\psi(s))^{\alpha-1}
\mathcal{D}^{\alpha,\psi}_{a+}f(s)ds.
\tag{3.7}\]
It follows that
\begin{align*}\label{eq3.8}
\bigg|f(t)-\sum_{k=0}^{n-1}\frac{f^{[k]}_{\psi}(a)}{k!}(\psi(t)-\psi(a))^{k}\bigg|
\leq &\quad\frac{1}{\Gamma(\alpha)}\int_{a}^{t}\psi^{'}(s)(\psi(t)-\psi(s))^{\alpha-1}
\big|\mathcal{D}^{\alpha,\psi}_{a+}f(s)\big|ds\\
\leq &\quad\frac{1}{\Gamma(\alpha)}\bigg(\int_{a}^{t}\psi^{'}(s)(\psi(t)-\psi(s))^{\alpha-1}ds\bigg)
\|\mathcal{D}^{\alpha,\psi}_{a+}f\|_{L_{\infty}}\\
=&\quad\frac{\|\mathcal{D}^{\alpha,\psi}_{a+}f\|_{L_{\infty}}}{\Gamma(\alpha+1)}(\psi(t)-\psi(a))^{\alpha}.\tag{3.8}
\end{align*}
Similarly, from the right $\psi$-Caputo fractional Taylor's formula we have
\begin{align*}\label{eq3.9}
\bigg|f(t)-\sum_{k=0}^{n-1}\frac{f^{[k]}_{\psi}(b)}{k!}(\psi(t)-\psi(b))^{k}\bigg|
\leq &\quad\frac{1}{\Gamma(\alpha)}\int_{t}^{b}\psi^{'}(s)(\psi(s)-\psi(t))^{\alpha-1}
\big|\mathcal{D}^{\alpha,\psi}_{b-}f(s)\big|ds\\
\leq &\quad\frac{\|\mathcal{D}^{\alpha,\psi}_{b-}f\|_{L_{\infty}}}{\Gamma(\alpha+1)}(\psi(b)-\psi(t))^{\alpha}.\tag{3.9}
\end{align*}
Set
\[\label{eq3.10}
\rho:=\mathrm{max}\bigg(\frac{\|\mathcal{D}^{\alpha,\psi}_{a+}f\|_{L_{\infty}}}{\Gamma(\alpha+1)},
\frac{\|\mathcal{D}^{\alpha,\psi}_{b-}f\|_{L_{\infty}}}{\Gamma(\alpha+1)}\bigg).\tag{3.10}
\]
Then from \eqref{eq3.8}, \eqref{eq3.9} and \eqref{eq3.10} we get
\[\label{eq3.11}
\bigg|f(t)-\sum_{k=0}^{n-1}\frac{f^{[k]}_{\psi}(a)}{k!}(\psi(t)-\psi(a))^{k}\bigg|
\leq \rho(\psi(t)-\psi(a))^{\alpha}\tag{3.11}
\]
and
\[\label{eq3.12}
\bigg|f(t)-\sum_{k=0}^{n-1}\frac{f^{[k]}_{\psi}(b)}{k!}(\psi(t)-\psi(b))^{k}\bigg|
\leq \rho(\psi(b)-\psi(t))^{\alpha}.\tag{3.12}
\]
It follows that
\begin{align*}\label{eq3.13}
\sum_{k=0}^{n-1}\frac{f^{[k]}_{\psi}(a)}{k!}(\psi(t)-\psi(a))^{k}-&\rho(\psi(t)-\psi(a))^{\alpha}\leq \quad f(t)\\
\leq &\quad\sum_{k=0}^{n-1}\frac{f^{[k]}_{\psi}(a)}{k!}(\psi(t)-\psi(a))^{k}+\rho(\psi(t)-\psi(a))^{\alpha}\tag{3.13}
\end{align*}
and
\begin{align*}\label{eq3.14}
\sum_{k=0}^{n-1}\frac{f^{[k]}_{\psi}(b)}{k!}(\psi(t)-\psi(b))^{k}-&\rho(\psi(b)-\psi(t))^{\alpha}\leq \quad f(t)\\
\leq &\quad\sum_{k=0}^{n-1}\frac{f^{[k]}_{\psi}(b)}{k!}(\psi(t)-\psi(b))^{k}+\rho(\psi(b)-\psi(t))^{\alpha}.\tag{3.14}
\end{align*}
Integrating \eqref{eq3.13} with respect to $\textit{t}$ from $\textit{a}$ to $\textit{s}$, we obtain
\begin{align*}\label{eq3.15}
\sum_{k=0}^{n-1}\frac{f^{[k]}_{\psi}(a)}{(k+1)!}&(\psi(s)-\psi(a))^{k+1}
-\frac{\rho}{(\alpha+1)}(\psi(s)-\psi(a))^{\alpha+1}\leq \int_{a}^{s}f(t)dt \\
\leq \quad &\sum_{k=0}^{n-1}\frac{f^{[k]}_{\psi}(a)}{(k+1)!}(\psi(s)-\psi(a))^{k+1}
+\frac{\rho}{(\alpha+1)}(\psi(s)-\psi(a))^{\alpha+1}.\tag{3.15}
\end{align*}
Integrating \eqref{eq3.14} with respect to $\textit{t}$ from $\textit{s}$ to $\textit{b}$, we obtain
\begin{align*}\label{eq3.16}
-\sum_{k=0}^{n-1}\frac{f^{[k]}_{\psi}(b)}{(k+1)!}&(\psi(s)-\psi(b))^{k+1}
-\frac{\rho}{(\alpha+1)}(\psi(b)-\psi(s))^{\alpha+1}\leq \int_{s}^{b}f(t)dt \\
\leq \quad &-\sum_{k=0}^{n-1}\frac{f^{[k]}_{\psi}(b)}{(k+1)!}(\psi(s)-\psi(b))^{k+1}
+\frac{\rho}{(\alpha+1)}(\psi(b)-\psi(s))^{\alpha+1}.\tag{3.16}
\end{align*}
Adding \eqref{eq3.15} and \eqref{eq3.16}, we get
\begin{align*}\label{eq3.17}
\sum_{k=0}^{n-1}\frac{1}{(k+1)!}&\big[f^{[k]}_{\psi}(a)(\psi(s)-\psi(a))^{k+1}-f^{[k]}_{\psi}(b)(\psi(s)-\psi(b))^{k+1}\big]\\
-&\frac{\rho}{(\alpha+1)}\big[(\psi(s)-\psi(a))^{\alpha+1}+(\psi(b)-\psi(s))^{\alpha+1}\big]
\leq\int_{a}^{b}f(t)dt\\
\leq \quad&\sum_{k=0}^{n-1}\frac{1}{(k+1)!}\big[f^{[k]}_{\psi}(a)(\psi(s)-\psi(a))^{k+1}-f^{[k]}_{\psi}(b)(\psi(s)-\psi(b))^{k+1}\big]\\
+&\frac{\rho}{(\alpha+1)}\big[(\psi(s)-\psi(a))^{\alpha+1}+(\psi(b)-\psi(s))^{\alpha+1}\big].\tag{3.17}
\end{align*}
From above it follows that
\begin{align*}\label{eq3.18}
\bigg|\int_{a}^{b}f(t)dt-\sum_{k=0}^{n-1}&\frac{1}{(k+1)!}\big[f^{[k]}_{\psi}(a)(\psi(s)-\psi(a))^{k+1}
+(-1)^{k}f^{[k]}_{\psi}(b)(\psi(s)-\psi(b))^{k+1}\big]\bigg|\\
\leq\quad &\frac{\rho}{(\alpha+1)}\big[(\psi(s)-\psi(a))^{\alpha+1}+(\psi(b)-\psi(s))^{\alpha+1}\big].\tag{3.18}
\end{align*}
(ii) Suppose ${\mathrm{u}(\psi(s))}=(\psi(s)-\psi(a))^{\alpha+1}+(\psi(b)-\psi(s))^{\alpha+1}$, for all $\psi(s)\in [\psi(a), \psi(b)]$.
Then
\begin{align*}
{\mathrm{u^{'}}(\psi(s))}=&(\alpha+1)\big[(\psi(s)-\psi(a))^{\alpha}-(\psi(b)-\psi(s))^{\alpha}\big]=0\\
&\Rightarrow \quad (\psi(s)-\psi(a))^{\alpha}=(\psi(b)-\psi(s))^{\alpha}\\
&\Rightarrow \quad \psi(s)-\psi(a)=\psi(b)-\psi(s)\\
&\Rightarrow \quad \psi(s)=\frac{\psi(a)+\psi(b)}{2}
\end{align*}
Therefore $\psi(s)=\frac{\psi(a)+\psi(b)}{2}$ is the critical point of $\mathrm{u}$.
If we put $\psi(s)=\frac{\psi(a)+\psi(b)}{2}$ in \eqref{eq3.1} we get the inequality \eqref{eq3.2}.\newline
(iii) When $f^{[k]}_{\psi}(a)=f^{[k]}_{\psi}(b)=0$, for $k=0,1,2,...,n-1,$ from \eqref{eq3.2} we get a sharp inequality as follows:
\begin{align*}\label{eq3.19}
\bigg|\int_{a}^{b}f(t)dt\bigg|
\leq \quad& \frac{\mathrm{max}\big[\|\mathcal{D}^{\alpha,\psi}_{a+}f\|_{L_{\infty}},\|\mathcal{D}^{\alpha,\psi}_{b-}f\|_{L_{\infty}}\big]}
{\Gamma(\alpha+2)}\frac{(\psi(b)-\psi(a))^{\alpha+1}}{2^{\alpha}}.\tag{3.19}
\end{align*}
(iv) Let $\psi(t_{i})=\psi(a)+i\bigg(\frac{\psi(b)-\psi(a)}{m}\bigg), i=0,1,2,...,m$.\newline
That is
\[\label{eq3.20}
\psi(t_{0})=\psi(a), \quad\psi(t_{1})=\psi(a)+\bigg(\frac{\psi(b)-\psi(a)}{m}\bigg),\,...\, ,\psi(t_{m})=\psi(b).\tag{3.20}
\]
Therefore
\[\label{eq3.21}
(\psi(t_{i})-\psi(a))=i\bigg(\frac{\psi(b)-\psi(a)}{m}\bigg),\quad(\psi(b)-\psi(t_{i}))=(m-i)\bigg(\frac{\psi(b)-\psi(a)}{m}\bigg).\tag{3.21}
\]
Therefore by \eqref{eq3.21} we can write
\begin{align*}\label{eq3.22}
(\psi(t_{i})-\psi(a))^{\alpha+1}+(\psi(b)-\psi(t_{i}))^{\alpha+1}=
\bigg(\frac{\psi(b)-\psi(a)}{m}\bigg)^{\alpha+1}\big(i^{\alpha+1}+(m-i)^{\alpha+1}\big).\tag{3.22}
\end{align*}
Also by using \eqref{eq3.21} for {\it k=0,1,2,...,n-1} and {\it i=0,1,2,...m} we have
\begin{align*}\label{eq3.23}
f^{[k]}(a)(\psi(t_{i})-\psi(a))^{k+1}+&(-1)^{k}f^{[k]}(b)(\psi(b)-\psi(t_{i}))^{k+1}\\
=&\bigg(\frac{\psi(b)-\psi(a)}{m}\bigg)^{k+1}[f^{[k]}(a)i^{k+1}+(-1)^{k}f^{[k]}(b)(m-i)^{k+1}].\tag{3.23}
\end{align*}
Hence from \eqref{eq3.18}, \eqref{eq3.22} and \eqref{eq3.23}, we get
\begin{align*}\label{eq3.24}
\bigg|\int_{a}^{b}f(t)dt-\sum_{k=0}^{n-1}&\frac{1}{(k+1)!}\big(\frac{\psi(b)-\psi(a)}{m}\big)^{k+1}\big[f^{[k]}_{\psi}(a)i^{k+1}
+(-1)^{k}f^{[k]}_{\psi}(b)(m-i)^{k+1}\big]\bigg|\\
\leq&\quad \frac{\rho}{(\alpha+1)}\bigg(\frac{\psi(b)-\psi(a)}{m}\bigg)^{\alpha+1}[i^{\alpha+1}+(m-i)^{\alpha+1}].\tag{3.24}
\end{align*}
From above we get required inequality.\newline
(v) If $f^{[k]}_{\psi}(a)=f^{[k]}_{\psi}(b)=0, k=1,2,...,n-1$, then from \eqref{eq3.24} we obtain
\begin{align*}\label{eq3.25}
\bigg|\int_{a}^{b}f(t)dt-&\bigg(\frac{\psi(b)-\psi(a)}{m}\bigg)\big[if(a)+(m-i)f(b)\big]\bigg|\\
\leq &\quad\frac{\rho}{(\alpha+1)}\bigg(\frac{\psi(b)-\psi(a)}{m}\bigg)^{\alpha+1}\big[i^{\alpha+1}+(m-i)^{\alpha+1}\big].\tag{3.25}
\end{align*}
(vi) For {\it i=1, m=2} from \eqref{eq3.25} we get
\begin{align*}
\bigg|\int_{a}^{b}f(t)dt-&\bigg(\frac{\psi(b)-\psi(a)}{2}\bigg)\big(f(a)+f(b)\big)\bigg|\\
\leq &\quad\frac{\rho}{(\alpha+1)}\frac{\big(\psi(b)-\psi(a)\big)^{\alpha+1}}{2^{\alpha}}.\qed
\end{align*}
{\textbf{Theorem 3.2}} Let $\alpha\geq1$, $\textit{f}\in \mathcal {\textit{A}C}^{n}([a, b])$ and $\psi \in \mathcal{C}^{n}([a, b])$ with $\psi$ is an increasing and $\psi^{'}(t)\neq0$, $\forall t\in [a, b]$. Suppose $\mathcal{D}^{\alpha, \psi}_{a+}f$, $\mathcal{D}^{\alpha, \psi}_{b-}f \in \mathcal{L}_{1}([a, b],\psi)$. Then for $a\leq s \leq b$ following inequalities hold:\newline
(i)
\begin{align*}\label{eq3.26}
&\bigg|\int_{a}^{b}f(t)dt-\sum_{k=0}^{n-1}\frac{1}{(k+1)!}
\big(f^{[k]}_{\psi}(a)(\psi(s)-\psi(a))^{k+1}+(-1)^{k}f^{[k]}_{\psi}(b)(\psi(b)-\psi(s))^{k+1}\big)\bigg| \\
\leq \quad& \frac{\mathrm{max}\big[\|\mathcal{D}^{\alpha,\psi}_{a+}f\|_{L_{1}([a, b],\psi)},
\|\mathcal{D}^{\alpha,\psi}_{b-}f\|_{L_{1}([a, b],\psi)}\big]}
{\Gamma(\alpha+2)}[(\psi(s)-\psi(a))^{\alpha+1}+(\psi(b)-\psi(s))^{\alpha+1}].\tag{3.26}
\end{align*}
(ii) The right hand side of the inequality \eqref{eq3.26} is minimized at $\psi(s)=\frac{\psi(a)+\psi(b)}{2}$, with the value $\frac{\mathrm{max}\big[\|\mathcal{D}^{\alpha,\psi}_{a+}f\|_{L_{1}([a, b],\psi)},
\|\mathcal{D}^{\alpha,\psi}_{b-}f\|_{L_{1}([a, b],\psi)}\big]}{\Gamma(\alpha+1)}\frac{(\psi(b)-\psi(a))^{\alpha}}{2^{\alpha-1}}$, that is
\begin{align*}\label{eq3.27}
\bigg|\int_{a}^{b}f(t)dt-&\sum_{k=0}^{n-1}\frac{1}{(k+1)!}\frac{(\psi(b)-\psi(a))^{k+1}}{2^{k+1}}
\big(f^{[k]}_{\psi}(a)+(-1)^{k}f^{[k]}_{\psi}(b)\big)\bigg|\\
\leq \quad& \frac{\mathrm{max}\big[\|\mathcal{D}^{\alpha,\psi}_{a+}f\|_{L_{1}([a, b],\psi)},
\|\mathcal{D}^{\alpha,\psi}_{b-}f\|_{L_{1}([a, b],\psi)}\big]}
{\Gamma(\alpha+1)}\frac{(\psi(b)-\psi(a))^{\alpha}}{2^{\alpha-1}}.\tag{3.27}
\end{align*}
(iii) When $f^{[k]}_{\psi}(a)=f^{[k]}_{\psi}(b)=0$, for $k=0,1,2,...,n-1,$ from \eqref{eq3.27} we get a sharp inequality as follows:
\begin{align*}\label{eq3.28}
\bigg|\int_{a}^{b}f(t)dt\bigg|
\leq \quad& \frac{\mathrm{max}\big[\|\mathcal{D}^{\alpha,\psi}_{a+}f\|_{L_{1}([a, b],\psi)},
\|\mathcal{D}^{\alpha,\psi}_{b-}f\|_{L_{1}([a, b],\psi)}\big]}
{\Gamma(\alpha+1)}\frac{(\psi(b)-\psi(a))^{\alpha}}{2^{\alpha-1}}.\tag{3.28}
\end{align*}
(iv) In general, for $i=0,1,2,...m\in\mathbb{N}$, we obtain
\begin{align*}\label{eq3.29}
\bigg|\int_{a}^{b}f(t)dt-\sum_{k=0}^{n-1}&\frac{1}{(k+1)!}\bigg(\frac{\psi(b)-\psi(a)}{m}\bigg)^{k+1}\big[i^{k+1}f^{[k]}_{\psi}(a)
+(-1)^{k}(m-i)^{k+1}f^{[k]}_{\psi}(b)\big]\bigg|\\
\leq&\quad \frac{\mathrm{max}\big[\|\mathcal{D}^{\alpha,\psi}_{a+}f\|_{L_{1}([a, b],\psi)},
\|\mathcal{D}^{\alpha,\psi}_{b-}f\|_{L_{1}([a, b],\psi)}\big]}
{\Gamma(\alpha+1)}\bigg(\frac{\psi(b)-\psi(a)}{m}\bigg)^{\alpha}\\
&\quad \big[i^{\alpha}+(m-i)^{\alpha}\big].\tag{3.29}
\end{align*}
(v) If $f^{[k]}_{\psi}(a)=f^{[k]}_{\psi}(b)=0, k=1,2,...,n-1$, then from \eqref{eq3.29} we obtain
\begin{align*}\label{eq3.30}
\bigg|\int_{a}^{b}f(t)dt-&\bigg(\frac{\psi(b)-\psi(a)}{m}\bigg)\big[if(a)+(m-i)f(b)\big]\bigg|\\
\leq &\quad\frac{\mathrm{max}\big[\|\mathcal{D}^{\alpha,\psi}_{a+}f\|_{L_{1}([a, b],\psi)},
\|\mathcal{D}^{\alpha,\psi}_{b-}f\|_{L_{1}([a, b],\psi)}\big]}
{\Gamma(\alpha+1)}\bigg(\frac{\psi(b)-\psi(a)}{m}\bigg)^{\alpha}\\
&\quad\big[i^{\alpha}+(m-i)^{\alpha}\big].\tag{3.30}
\end{align*}
(vi) For {\it i=1, m=2} from \eqref{eq3.30} we get
\begin{align*}\label{eq3.31}
\bigg|\int_{a}^{b}f(t)dt-&\bigg(\frac{\psi(b)-\psi(a)}{2}\bigg)\big(f(a)+f(b)\big)\bigg|\\
\leq &\quad \frac{\mathrm{max}\big[\|\mathcal{D}^{\alpha,\psi}_{a+}f\|_{L_{1}([a, b],\psi)},
\|\mathcal{D}^{\alpha,\psi}_{b-}f\|_{L_{1}([a, b],\psi)}\big]}
{\Gamma(\alpha+1)}\frac{\big(\psi(b)-\psi(a)\big)^{\alpha}}{2^{\alpha-1}}.\tag{3.31}
\end{align*}
\textbf{Proof} From the left $\psi$-Caputo fractional Taylor's formula we have
\[\label{eq3.32}
f(t)-\sum_{k=0}^{n-1}\frac{f^{[k]}_{\psi}(a)}{k!}(\psi(t)-\psi(a))^{k}
=\frac{1}{\Gamma(\alpha)}\int_{a}^{t}\psi^{'}(s)(\psi(t)-\psi(s))^{\alpha-1}
\mathcal{D}^{\alpha,\psi}_{a+}f(s)ds.\tag{3.32}
\]
Then
\begin{align*}\label{eq3.33}
\bigg|f(t)-\sum_{k=0}^{n-1}\frac{f^{[k]}_{\psi}(a)}{k!}(\psi(t)-\psi(a))^{k}\bigg|
\leq &\quad\frac{1}{\Gamma(\alpha)}\int_{a}^{t}\psi^{'}(s)(\psi(t)-\psi(s))^{\alpha-1}
\big|\mathcal{D}^{\alpha,\psi}_{a+}f(s)\big|ds.\tag{3.33}\\
\end{align*}
Here $\alpha\geq1$ and $\psi$ is increasing for $a\leq s \leq t$, we have
\begin{align*}\label{eq3.34}
\bigg|f(t)-\sum_{k=0}^{n-1}\frac{f^{[k]}_{\psi}(a)}{k!}(\psi(t)-\psi(a))^{k}\bigg|
\leq &\quad\frac{(\psi(t)-\psi(a))^{\alpha-1}}{\Gamma(\alpha)}\int_{a}^{t}\|\mathcal{D}^{\alpha,\psi}_{a+}f\|_{L_{1}}d\psi(s)\\
=&\quad\frac{\|\mathcal{D}^{\alpha,\psi}_{a+}f\|_{L_{1}([a, b],\psi)}}{\Gamma(\alpha)}(\psi(t)-\psi(a))^{\alpha-1}.\tag{3.34}
\end{align*}
Similarly, from the right $\psi$-Caputo fractional Taylor's formula we obtain
\begin{align*}\label{eq3.35}
\bigg|f(t)-\sum_{k=0}^{n-1}\frac{f^{[k]}_{\psi}(b)}{k!}(\psi(t)-\psi(b))^{k}\bigg|
\leq&\quad\frac{\|\mathcal{D}^{\alpha,\psi}_{b-}f\|_{L_{1}([a, b],\psi)}}{\Gamma(\alpha)}(\psi(b)-\psi(t))^{\alpha-1}.\tag{3.35}
\end{align*}
Set
\[\label{eq3.36}
\eta:=\mathrm{max}\bigg(\frac{\|\mathcal{D}^{\alpha,\psi}_{a+}f\|_{L_{1}([a, b],\psi)}}{\Gamma(\alpha)},
\frac{\|\mathcal{D}^{\alpha,\psi}_{b-}f\|_{L_{1}([a, b],\psi)}}{\Gamma(\alpha)}\bigg).\tag{3.36}
\]
Hence from \eqref{eq3.34}, \eqref{eq3.35} and \eqref{eq3.36} we have
\[\label{eq3.37}
\bigg|f(t)-\sum_{k=0}^{n-1}\frac{f^{[k]}_{\psi}(a)}{k!}(\psi(t)-\psi(a))^{k}\bigg|
\leq \eta(\psi(t)-\psi(a))^{\alpha-1}\tag{3.37}
\]
and
\[\label{eq3.38}
\bigg|f(t)-\sum_{k=0}^{n-1}\frac{f^{[k]}_{\psi}(b)}{k!}(\psi(t)-\psi(b))^{k}\bigg|
\leq \eta(\psi(b)-\psi(t))^{\alpha-1}.\tag{3.38}
\]
By using similar arguments as in the proof of theorem (3.1), we can write
\begin{align*}
\bigg|\int_{a}^{b}f(t)dt-\sum_{k=0}^{n-1}&\frac{1}{(k+1)!}\big[f^{[k]}_{\psi}(a)(\psi(s)-\psi(a))^{k+1}
+(-1)^{k}f^{[k]}_{\psi}(b)(\psi(s)-\psi(b))^{k+1}\big]\bigg|\\
\leq\quad &\frac{\eta}{\alpha}\big[(\psi(s)-\psi(a))^{\alpha}+(\psi(b)-\psi(s))^{\alpha}\big].
\end{align*}
The rest of the proof similar to the proof of theorem (3.1).$\qed$\\
{\textbf{Theorem 3.3}} Let $ p, q > 1$ with $\frac{1}{p}+\frac{1}{q}=1$ and $\alpha>\frac{1}{q}>0$. Let $\textit{f}\in \mathcal {\textit{A}C}^{n}([a, b])$ and $\psi \in \mathcal{C}^{n}([a, b])$ with $\psi$ is an increasing and $\psi^{'}(t)\neq0$, $\forall t\in [a, b]$. Suppose $\mathcal{D}^{\alpha, \psi}_{a+}f$, $\mathcal{D}^{\alpha, \psi}_{b-}f \in \mathcal{L}_{q}([a, b],\psi)$. Then for $a\leq s \leq b$ following inequalities hold:\newline
(i)
\begin{align*}\label{eq3.39}
&\bigg|\int_{a}^{b}f(t)dt-\sum_{k=0}^{n-1}\frac{1}{(k+1)!}
\big(f^{[k]}_{\psi}(a)(\psi(s)-\psi(a))^{k+1}+(-1)^{k}f^{[k]}_{\psi}(b)(\psi(b)-\psi(s))^{k+1}\big)\bigg| \\
\leq \quad& \frac{\mathrm{max}\big[\|\mathcal{D}^{\alpha,\psi}_{a+}f\|_{L_{q}([a, b],\psi)},\|
\mathcal{D}^{\alpha,\psi}_{b-}f\|_{L_{q}([a, b],\psi)}\big]}
{\Gamma(\alpha)(\alpha+\frac{1}{p})(p(\alpha-1)+1)^{\frac{1}{p}}}
[(\psi(s)-\psi(a))^{\alpha+\frac{1}{p}}+(\psi(b)-\psi(s))^{\alpha+\frac{1}{p}}].\tag{3.39}
\end{align*}
(ii) The right hand side of the inequality \eqref{eq3.39} is minimized at $\psi(s)=\frac{\psi(a)+\psi(b)}{2}$, with the value $\frac{\mathrm{max}\big[\|\mathcal{D}^{\alpha,\psi}_{a+}f\|_{L_{q}([a, b],\psi)},
\|\mathcal{D}^{\alpha,\psi}_{b-}f\|_{L_{q}([a, b],\psi)}\big]}
{\Gamma(\alpha)(\alpha+\frac{1}{p})(p(\alpha-1)+1)^{\frac{1}{p}}}
\frac{(\psi(b)-\psi(a))^{\alpha+\frac{1}{p}}}{2^{\alpha-\frac{1}{q}}}$, that is
\begin{align*}\label{eq3.40}
\bigg|\int_{a}^{b}f(t)dt-&\sum_{k=0}^{n-1}\frac{1}{(k+1)!}\frac{(\psi(b)-\psi(a))^{k+1}}{2^{k+1}}
\big(f^{[k]}_{\psi}(a)+(-1)^{k}f^{[k]}_{\psi}(b)\big)\bigg|\\
\leq \quad& \frac{\mathrm{max}\big[\|\mathcal{D}^{\alpha,\psi}_{a+}f\|_{L_{q}([a, b],\psi)},
\|\mathcal{D}^{\alpha,\psi}_{b-}f\|_{L_{q}([a, b],\psi)}\big]}
{\Gamma(\alpha)(\alpha+\frac{1}{p})(p(\alpha-1)+1)^{\frac{1}{p}}}
\frac{(\psi(b)-\psi(a))^{\alpha+\frac{1}{p}}}{2^{\alpha-\frac{1}{q}}}.\tag{3.40}
\end{align*}
(iii) When $f^{[k]}_{\psi}(a)=f^{[k]}_{\psi}(b)=0$, for $k=0,1,2,...,n-1,$ from \eqref{eq3.40} we get a sharp inequality as follows:
\begin{align*}\label{eq3.41}
\bigg|\int_{a}^{b}f(t)dt\bigg|
\leq \quad& \frac{\mathrm{max}\big[\|\mathcal{D}^{\alpha,\psi}_{a+}f\|_{L_{q}([a, b],\psi)},
\|\mathcal{D}^{\alpha,\psi}_{b-}f\|_{L_{q}([a, b],\psi)}\big]}
{\Gamma(\alpha)(\alpha+\frac{1}{p})(p(\alpha-1)+1)^{\frac{1}{p}}}
\frac{(\psi(b)-\psi(a))^{\alpha+\frac{1}{p}}}{2^{\alpha-\frac{1}{q}}}.\tag{3.41}
\end{align*}
(iv) In general, for $i=0,1,2,...m\in\mathbb{N}$, we obtain
\begin{align*}\label{eq3.42}
\bigg|\int_{a}^{b}f(t)dt-\sum_{k=0}^{n-1}&\frac{1}{(k+1)!}\bigg(\frac{\psi(b)-\psi(a)}{m}\bigg)^{k+1}\big[i^{k+1}f^{[k]}_{\psi}(a)
+(-1)^{k}(m-i)^{k+1}f^{[k]}_{\psi}(b)\big]\bigg|\\
\leq&\quad \frac{\mathrm{max}\big[\|\mathcal{D}^{\alpha,\psi}_{a+}f\|_{L_{q}([a, b],\psi)},
\|\mathcal{D}^{\alpha,\psi}_{b-}f\|_{L_{q}([a, b],\psi)}\big]}
{\Gamma(\alpha)(\alpha+\frac{1}{p})(p(\alpha-1)+1)^{\frac{1}{p}}}
\bigg(\frac{\psi(b)-\psi(a)}{m}\bigg)^{\alpha+\frac{1}{p}}\\
&\quad[i^{\alpha+\frac{1}{p}}+(m-i)^{\alpha+\frac{1}{p}}].\tag{3.42}
\end{align*}
(v) If $f^{[k]}_{\psi}(a)=f^{[k]}_{\psi}(b)=0, k=1,2,...,n-1$, then from \eqref{eq3.42} we obtain
\begin{align*}\label{eq3.43}
\bigg|\int_{a}^{b}f(t)dt-&\bigg(\frac{\psi(b)-\psi(a)}{m}\bigg)\big[if(a)+(m-i)f(b)\big]\bigg|\\
\leq &\quad\frac{\mathrm{max}\big[\|\mathcal{D}^{\alpha,\psi}_{a+}f\|_{L_{q}([a, b],\psi)},
\|\mathcal{D}^{\alpha,\psi}_{b-}f\|_{L_{q}([a, b],\psi)}\big]}
{\Gamma(\alpha)(\alpha+\frac{1}{p})(p(\alpha-1)+1)^{\frac{1}{p}}}
\bigg(\frac{\psi(b)-\psi(a)}{m}\bigg)^{\alpha+\frac{1}{p}}\\
&\quad\big[i^{\alpha+\frac{1}{p}}+(m-i)^{\alpha+\frac{1}{p}}\big].\tag{3.43}
\end{align*}
(vi) For {\it i=1, m=2} from \eqref{eq3.43} we get
\begin{align*}\label{eq3.44}
\bigg|\int_{a}^{b}f(t)dt-&\bigg(\frac{\psi(b)-\psi(a)}{2}\bigg)\big(f(a)+f(b)\big)\bigg|\\
\leq &\quad \frac{\mathrm{max}\big[\|\mathcal{D}^{\alpha,\psi}_{a+}f\|_{L_{q}([a, b],\psi)},
\|\mathcal{D}^{\alpha,\psi}_{b-}f\|_{L_{q}([a, b],\psi)}\big]}
{\Gamma(\alpha)(\alpha+\frac{1}{p})(p(\alpha-1)+1)^{\frac{1}{p}}}
\frac{\big(\psi(b)-\psi(a)\big)^{\alpha+\frac{1}{p}}}{2^{\alpha-\frac{1}{q}}}.\tag{3.44}
\end{align*}

\textbf{Proof} From left sided $\psi$-Caputo fractional Taylor's formula we have
\begin{align*}
\bigg|f(t)-\sum_{k=0}^{n-1}\frac{f^{[k]}_{\psi}(a)}{k!}(\psi(t)-\psi(a))^{k}\bigg|
\leq &\quad\frac{1}{\Gamma(\alpha)}\int_{a}^{t}\psi^{'}(s)(\psi(t)-\psi(s))^{\alpha-1}
\big|\mathcal{D}^{\alpha,\psi}_{a+}f(s)\big|ds.
\end{align*}
By using H\"{o}lder's inequality, we obtain
\begin{align*}\label{eq3.45}
\bigg|f(t)-\sum_{k=0}^{n-1}\frac{f^{[k]}_{\psi}(a)}{k!}(\psi(t)-\psi(a))^{k}\bigg|
\leq &\quad\frac{1}{\Gamma(\alpha)}\bigg(\int_{a}^{t}(\psi(t)-\psi(s))^{p(\alpha-1)}d\psi(s)\bigg)^{\frac{1}{p}}\\
&\bigg(\int_{a}^{t}\big|\mathcal{D}^{\alpha,\psi}_{a+}f(s)d\psi(s)\big|^{q}\bigg)^{\frac{1}{q}}\\
\leq &\quad\frac{1}{\Gamma(\alpha)}\frac{(\psi(t)-\psi(a))^{\frac{p(\alpha-1)+1}{p}}}{(p(\alpha-1)+1)^{\frac{1}{p}}}
\|\mathcal{D}^{\alpha,\psi}_{a+}f\|_{L_{q}([a, b],\psi)}\\
=&\quad\frac{1}{\Gamma(\alpha)}\frac{\|\mathcal{D}^{\alpha,\psi}_{a+}f\|_{L_{q}([a, b],\psi)}}{(p(\alpha-1)+1)^{\frac{1}{p}}}
(\psi(t)-\psi(a))^{\alpha-\frac{1}{q}}.\tag{3.45}
\end{align*}
Similarly, from right sided $\psi$-Caputo fractional Taylor's formula we have
\begin{align*}\label{eq3.46}
\bigg|f(t)-\sum_{k=0}^{n-1}\frac{f^{[k]}_{\psi}(b)}{k!}(\psi(t)-\psi(b))^{k}\bigg|
\leq &\quad\frac{1}{\Gamma(\alpha)}\frac{\|\mathcal{D}^{\alpha,\psi}_{b-}f\|_{[L_{q},\psi]}}{(p(\alpha-1)+1)^{\frac{1}{p}}}
(\psi(b)-\psi(t))^{\alpha-\frac{1}{q}}.\tag{3.46}
\end{align*}
Set
\[\label{eq3.47}
\delta:=\mathrm{max}\bigg(\frac{\|\mathcal{D}^{\alpha,\psi}_{a+}f\|_{L_{q}([a, b],\psi)}}{\Gamma(\alpha)},
\frac{\|\mathcal{D}^{\alpha,\psi}_{b-}f\|_{L_{q}([a, b],\psi)}}{\Gamma(\alpha)}\bigg).\tag{3.47}
\]
From \eqref{eq3.45}, \eqref{eq3.46} and \eqref{eq3.47} it follows that
\[
\bigg|f(t)-\sum_{k=0}^{n-1}\frac{f^{[k]}_{\psi}(a)}{k!}(\psi(t)-\psi(a))^{k}\bigg|\leq \delta(\psi(t)-\psi(a))^{\alpha-\frac{1}{q}}
\]
and
\[
\bigg|f(t)-\sum_{k=0}^{n-1}\frac{f^{[k]}_{\psi}(b)}{k!}(\psi(t)-\psi(b))^{k}\bigg|
\leq \delta(\psi(b)-\psi(t))^{\alpha-\frac{1}{q}}.
\]
The remaining proof follows by using the similar arguments as in the proof of theorem (3.1).$\qed$\\

\end{document}